\documentclass[12pt, twoside, eqno]{article}
\usepackage{latexsym}
\usepackage{amssymb}
\usepackage{amsfonts}
\begin{document}
\begin{center}
{\bf \large On fuzzy prime and fuzzy semiprime ideals of 
$\le$-hypergroupoids}\medskip

\noindent{\bf Niovi Kehayopulu}\medskip

{\small Department of Mathematics, University of Athens,
15784 Panepistimiopolis,\\Athens, Greece\\email: 
nkehayop@math.uoa.gr}

\end{center}
\bigskip

{\small
\noindent{\bf Abstract.} We deal with an hypergroupoid endowed with a 
relation denoted by ``$\le$", we call it $\le$--hypergroupoid. We 
prove that a nonempty subset $A$ of a $\le$--hypergroupoid $H$ is a 
prime (resp. semiprime) ideal of $H$ if and only if its 
characteristic function $f_A$ is a fuzzy prime (resp. fuzzy 
semiprime) ideal of $H$.\medskip

}
\section{Introduction and prerequisites}
A characterization of prime and semiprime ideals of semigroups in 
terms of fuzzy subsets has been considered in [1], and similar 
characterizations hold for ordered groupoids in general. Fuzzy sets 
in ordered groupoids have been first considered in [2]. In the 
present paper we examine the results in [1] in case of an 
hypergroupoid $H$ endowed with a relation denoted by ``$\le$" (not an 
ordered relation, as so no compatible with the multiplication of $H$ 
in general). As a consequence, our results hold for ordered 
hypergroupoids as well. An {\it hypergroupoid} is a nonempty set $H$ 
with an hyperoperation $$\circ : H\times H \rightarrow {\cal P}^*(H) 
\mid (a,b) \rightarrow a\circ b$$on $H$ and an operation $$* : {\cal 
P}^*(H)\times {\cal P}^*(H) \rightarrow {\cal P}^*(H) \mid (A,B) 
\rightarrow A*B$$ on ${\cal P}^*(H)$ (induced by the operation of 
$H$) such that $$A*B=\bigcup\limits_{(a,b) \in\,A\times B} {(a\circ 
b)}$$ for every $A,B\in {\cal P}^*(H)$. An hypergroupoid can be also 
denoted by $(H,\circ)$ as the operation ``$*$" depends on ``$\circ$". 
A nonempty subset $A$ of an hypergroupoid $H$ is called a {\it left} 
(resp. {\it right}) {\it ideal} of $H$ if $H*A\subseteq A$ (resp. 
$A*H\subseteq A$). It is called an ideal of $H$ if it is both a left 
and a right ideal of $H$. If $H$ is an hypergroupoid then, for every 
$x,y\in H$, we have $\{x\}*\{y\}=x\circ y.$ The following 
proposition, though clear, plays an essential role in the theory of 
hypergroupoids.\medskip

\noindent{\bf Proposition 1.} {\it Let $(H,\circ)$ be an 
hypergroupoid, $x\in H$ and $A,B\in {\cal P}^*(H)$. Then we have the 
following:

$(1)$ $x\in A*B$ $\Longleftrightarrow$ $x\in a\circ b$ for some $a\in 
A$, $b\in B$.

$(2)$ If $a\in A$ and $b\in B$, then $a\circ b\subseteq A*B$.} 
\medskip

\noindent{\bf Proposition 2.} {\it Let $(H,\circ)$ be an 
hypergroupoid. If $A$ is a left (resp. right) ideal of H, then for 
every $h\in H$ and every $a\in A$, we have $h\circ a\subseteq A$ 
(resp. $a\circ h\subseteq A$). ``Conversely", if $A$ is a nonempty 
subset of $H$ such that $h\circ a\subseteq A$ (resp. $a\circ 
h\subseteq A$) for every $h\in H$ and every $a\in A$, then the set 
$A$ is a left (resp. right) ideal of H}.
\section{Main results}
\noindent{\bf Definition 3.} By a {\it $\le$--hypergroupoid} we mean 
an hypergroupoid $H$ endowed with a relation denoted by ``$\le$". 
\medskip

\noindent{\bf Definition 4.} Let $H$ be a $\le$--hypergroupoid. A 
nonempty subset $A$ of $H$ is called a {\it left} (resp. {\it right}) 
{\it ideal} of $H$ if

$(1)$ $H*A\subseteq A$ (resp. $A*H\subseteq A)$ and

$(2)$ if $a\in A$ and $H\ni b\le a$, then $b\in A$.\\A subset of $H$ 
which is both a left ideal and a right ideal of $H$ is called an {\it 
ideal} of $H$. A nonempty subset $A$ of $H$ is called a {\it 
subgroupoid} of $H$ if $A*A\subseteq A$.\\ Clearly, every left ideal, 
right ideal or ideal of $H$ is a subgroupoid of $H$.\medskip

\noindent{\bf Definition 5.} Let $H$ be an hypergroupoid (or a 
$\le$--hypergroupoid). A nonempty subset $I$ of $H$ is called a {\it 
prime subset} of $H$ if

(1) $a,b\in H \mbox { such that } a\circ b\subseteq I \mbox { implies 
} a\in I \mbox { or } b\in I$ and

(2) if $a,b\in H$, then $a\circ b\subseteq I$ or $(a\circ b)\cap 
I=\emptyset$.\\The following are equivalent:

$(1)$ $a,b\in H$, $a\circ b\subseteq I$ $\Longrightarrow$ $a\in I$ or 
$b\in I$.

$(2)$ $\emptyset\not=A,B\subseteq H$, $A*B\subseteq I$ 
$\Longrightarrow$ $A\subseteq I$ or $B\subseteq I$.\\Indeed: 
$(1)\Longrightarrow (2)$. Let $A,B\in {\cal P}^*(H)$, $A*B\subseteq 
I$ and $A\nsubseteq I$. Let $a\in A$ such that $a\notin I$ and $b\in 
B$ $(B\not=\emptyset)$. We have $a\circ b\subseteq A*B\subseteq I$. 
Then, by (1), $a\in I$ or $b\in I$.\\$(2)\Longrightarrow (1)$. Let 
$a,b\in H$, $a\circ b\subseteq I$. Then $\{a\}*\{b\}=a\circ 
b\subseteq I$. By (2), we have $\{a\}\subseteq I$ or $\{b\}\subseteq 
I$, so $a\in I$ or $b\in I$.\smallskip

By a prime ideal of $H$ we clearly mean an ideal of $H$ which is at 
the same time a prime subset of $H$.

Following Zadeh, any mapping $f : H\rightarrow [0,1]$ of a 
$\le$--hypergroupoid $H$ into the closed interval $[0,1]$ of real 
numbers is called a {\it fuzzy subset} of $H$ or a ({\it fuzzy set} 
in $H$) and $f_A$ (: the characteristic function) is the mapping$$f_A 
: H \rightarrow \{0,1\} \mid x \rightarrow f_A (x)=\left\{ 
\begin{array}{l}
1\,\,\,\,\,$if$\,\,\,\,x \in A\\
0\,\,\,\,$if$\,\,\,\,x \notin A.
\end{array} \right.$$
{\bf Definition 6.} Let $H$ be a $\le$--hypergroupoid. A fuzzy subset 
$f$ of $H$ is called a {\it fuzzy left ideal} of $H$ if

$(1)$ $x\le y \Rightarrow f(x)\ge f(y)$ and

$(2)$ if $f(x\circ y)\ge f(y)$ for all $x,y\in H$.\\With the property 
(2) we mean the following:

(2) if $x,y\in H$ and $u\in x\circ y$, then $f(u)\ge f(y)$.

A fuzzy subset $f$ of $H$ is called a {\it fuzzy right ideal} of $H$ 
if

$(1)$ $x\le y \Rightarrow f(x)\ge f(y)$ and

$(2)$ if $f(x\circ y)\ge f(x)$ for all $x,y\in H$.\\With the property 
(2) we mean the following:

(2) if $x,y\in H$ and $u\in x\circ y$, then $f(u)\ge f(x)$.

A fuzzy subset of $H$ is called a {\it fuzzy ideal} of $H$ if it is 
both a fuzzy left and a fuzzy right ideal of $H$. As one can easily 
see, a fuzzy subset $f$ of $H$ is a fuzzy ideal of $H$ if and only 
if$$f(x\circ y)\ge \max\{f(x),f(y)\} \mbox { for all } x,y\in H$$in 
the sense that$$x,y\in H \mbox { and } u\in x\circ y \mbox { implies 
} f(u)\ge \max\{f(x),f(y)\}.$${\bf Proposition 7.} {\it Let H be a 
$\le$--hypergroupoid. If $A$ is a left (resp. right) ideal of $H$, 
then the characteristic function $f_A$ is a fuzzy left (resp. fuzzy 
right) ideal of $H$. ``Conversely", if $A$ is a nonempty subset of 
$H$ such that $f_A$ is a fuzzy left (resp. fuzzy right) ideal of H, 
then the set A is a left (resp. right) ideal of H}. \medskip

\noindent{\bf Proposition 8.} {\it Let H be an $\le$--hypergroupoid. 
If $A$ is an ideal of $H$, then $f_A$ is a fuzzy ideal of $H$. 
``Conversely", if $A$ is a nonempty subset of $H$ such that $f_A$ is 
a fuzzy ideal of H, then the set A is an ideal of H}.
\medskip

\noindent{\bf Definition 9.} Let $H$ be an hypergroupoid (or a 
$\le$--hypergroupoid). A fuzzy subset $f$ of $H$ is called {\it fuzzy 
prime subset} of $H$ if
$$f(x\circ y)\le\max\{f(x),f(y)\} \mbox { for all } x,y\in H$$that 
is, if $x,y\in H$ and $u\in x\circ y$, then $f(u)\le 
\max\{f(x),f(y)\}$.\smallskip

By a fuzzy prime ideal of $H$ we mean a fuzzy ideal of $H$ which is 
at the same time a fuzzy prime subset of $H$. So a fuzzy subset $f$ 
of a $\le$--hypergroupoid $H$ is a fuzzy prime ideal of $H$ if and 
only if the following assertions are satisfied:

(1) $x\le y$ implies $f(x)\ge f(y)$ and

(2) $f(x\circ y)=\max\{f(x),f(y)\} \mbox { for all } x,y\in H$\\that 
is, if $x,y\in H$ and $u\in x\circ y$, then 
$f(u)=\max\{f(x),f(y)\}$.\medskip

\noindent{\bf Proposition 10.} {\it Let H be an $\le$--hypergroupoid. 
If A is a prime ideal of H, then $f_A$ is a fuzzy prime ideal of $H$. 
``Conversely", if $A$ is a nonempty subset of $H$ such that $f_A$ is 
a fuzzy prime ideal of H, then A is a prime ideal of H}.\medskip

\noindent{\bf Proof.} $\Longrightarrow$. Since $A$ is an ideal of 
$H$, $f_A$ is a fuzzy ideal of $H$. Let $x,y\in H$ and $u\in x\circ 
y$. Then $f_A(u)=\max\{f_A(x),f_A(y)\}$. Indeed: Let $x\circ 
y\subseteq A$. Since $A$ is a prime ideal of $H$, we have $x\in A$ or 
$y\in A$. Then $f_A(x)=1$ or $f_A(y)=1$, and 
$\max\{f_A(x),f_A(y)\}=1$. Since $u\in x\circ y\subseteq A$, we have 
$u\in A$. Then $f_A(u)=1$, so $f_A(u)=\max\{f_A(x),f_A(y)\}$. Let 
$x\circ y\nsubseteq A$. Since $A$ is a prime ideal of $H$, we have 
$(x\circ y)\cap A=\emptyset$. Since $u\in x\circ y$, we have $u\notin 
A$, so $f_A(u)=0$. Since $x\circ y\nsubseteq A$ and $A$ is an ideal 
of $H$, we have $x\notin A$ and $y\notin A$ (since $x\in A$ implies 
$x\circ y\subseteq A*H\subseteq A$ and $y\in A$ implies $x\circ 
y\subseteq H*A\subseteq A$ which is impossible). Then we have 
$f_A(x)=0=f_A(y)$, and $f_A(u)=\max\{f_A(x),f_A(y)\}$.\\
$\Longleftarrow$. Let $f_A$ be a fuzzy prime ideal of $H$. Since 
$f_A$ is a fuzzy ideal of $H$, $A$ is an ideal of $H$. Let $x,y\in H$ 
such that $x\circ y\subseteq A$. Suppose $x\notin A$ and $y\notin A$. 
Then $f_A(x)=0=f_A(y)$. Take an element $u\in x\circ y$ $(x\circ 
y\not=\emptyset)$. Since $u\in A$, we have $f_A(u)=1$, so
$f_A(u)\not=\max\{f_A(x),f_A(y)\}$ which is impossible. Thus we have 
$x\in A$ or $y\in A$. Let now $x,y\in H$ such that $x\circ 
y\nsubseteq A$. Then $(x\circ y)\cap A=\emptyset$. Indeed: Let $u\in 
(x\circ y)\cap A$. Since $u\in x\circ y$, by hypothesis, we have 
$f_A(u)=\max\{f_A(x),f_A(y)\}$. Since $u\in A$, we have $f_A(u)=1$. 
Then $ f_A(x)=1$ or $f_A(y)=1$, so $x\in A$ or $y\in A$. If $x\in A$, 
then $x\circ y\subseteq A*H\subseteq A$ (since $A$ is an ideal of 
$H$), which is impossible. If $y\in A$, then $x\circ y\subseteq 
H*A\subseteq A$ which again is impossible. Hence we have $(x\circ 
y)\cap A=\emptyset$. $\hfill\Box$\medskip

\noindent{\bf Definition 11.} Let $H$ be an hypergroupoid (or a 
$\le$--hypergroupoid). A nonempty subset $I$ of $H$ is called {\it 
semiprime subset} of $H$ if

(1) if $a\in H \mbox { such that } a\circ a\subseteq I$, then $a\in 
I$ and

(2) if $a\in H$, then $a\circ a\subseteq I$ or $(a\circ a)\cap 
I=\emptyset$.\\The following are equivalent:

$(1)$ if $a\in H$ such that $a\circ a\subseteq I$, then $a\in I$.

$(2)$ if $A$ is a nonempty subset of $H$ such that $A*A\subseteq I$, 
then $A\subseteq I$.

By a semiprime ideal of $H$ we clearly mean an ideal of $H$ which is 
at the same time a semiprime subset of $H$.\medskip

\noindent{\bf Definition 12.} Let $H$ be an hypergroupoid (or a 
$\le$--hypergroupoid). A fuzzy subset $f$ of $H$ is called {\it fuzzy 
semiprime subset} of $H$ if$$f(x)\ge f(x\circ x) \mbox { for every } 
x\in H$$that is, if $x\in H$ and $u\in x\circ x$, then $f(x)\ge 
f(u)$.

By a fuzzy semiprime ideal of $H$ we cleary mean a fuzzy ideal of $H$ 
which is at the same time a semiprime fuzzy subset of $H$.\medskip

\noindent{\bf Remark 13.} If $f$ is a fuzzy ideal of $H$ and $a\in 
H$, then $f(a\circ a)\ge \max\{f(a),f(a)\}=f(a)$. Hence: If $f$ is a 
fuzzy semiprime ideal of $H$, then $f(a\circ a)=f(a)$ for every $a\in 
H$. If $f$ is a fuzzy prime ideal of $H$ and $a\in H$, then $f(a\circ 
a)=\max\{f(a),f(a)\}=f(a)$, so $f(a)\ge f(a\circ a)$, and $f$ is 
fuzzy semiprime ideal.

A fuzzy subset $f$ of $H$ is a fuzzy semiprime ideal of $H$ if and 
only if the following assertions are satisfied:

(1) $x\le y$ implies $f(x)\ge f(y)$ and

(2) if $f(x\circ x)=f(x)$ for every $x\in H$\\that is, if $x\in H$ 
and $u\in x\circ x$, then $f(u)=f(x)$.
\medskip

\noindent{\bf Proposition 14.} {\it Let H be a $\le$--hypergroupoid. 
If A is a semiprime ideal of H, then $f_A$ is a fuzzy semiprime ideal 
of $H$. ``Conversely", if $A$ is a nonempty subset of $H$ such that 
$f_A$ is a fuzzy semiprime ideal of H, then A is a prime ideal of 
H}.\medskip

\noindent{\bf Proof.} $\Longrightarrow$. Let $A$ be a semiprime ideal 
of $H$. Since $A$ is an ideal of $H$, $f_A$ is a fuzzy ideal of $H$. 
Let $x\in H$ and $u\in x\circ x$. Then $f_A(u)=f_A(x)$. Indeed: Let 
$x\circ x\nsubseteq A$. Since $A$ is a semiprime subset of $H$, we 
have $(x\circ x)\cap A=\emptyset$, so $u\notin A$, and $f_A(u)=0$. On 
the other hand, since $x\circ x\not\subseteq A$ and $A$ is an ideal 
of $H$, we have $x\not\in A$, then $f_A(x)=0$, so $f_A(u)=f_A(x)$. 
Let $x\circ x\subseteq A$. Then $u\in A$, so $f_A(u)=1$. On the other 
hand, since $A$ is a semiprime subset of $H$ and $x\circ x\subseteq 
A$, we have $x\in A$, so $f_A(x)=1$. Then $f_A(u)= f_A(x)$.\\
\noindent $\Longleftarrow$. Let $f_A$ be a fuzzy semiprime ideal of 
$H$. Since $f_A$ is a fuzzy ideal of $H$, the set $A$ is an ideal of 
$H$. Let $x\in H$ such that $x\circ x\subseteq A$. Then $x\in A$. 
Indeed: Let $x\notin A$. Then $f_A(x)=0$. Take an element $u\in 
x\circ x$ $(x\circ x\not=\emptyset$). Since $u\in A$, we have 
$f_A(u)=1$. Since $f_A$ is a semiprime ideal of $H$, we have 
$f_A(u)=f_A(x)$ which is impossible. Thus we have $x\in A$. Let $x\in 
H$ such that $x\circ x\nsubseteq A$. Then $(x\circ x)\cap 
A=\emptyset$. Indeed: Let $u\in (x\circ x)\cap A$. Since $u\in x\circ 
x$, by hypothesis, we have $f_A(u)=f_A(x)$. Since $u\in A$, we have 
$f_A(u)=1$, then $f_A(x)=1$, and $x\in A$. Then $x\circ x\subseteq 
A*A\subseteq A$ (since $A$ is a subgroupoid of $H$), which is 
impossible. Thus we have $(x\circ x)\cap A=\emptyset$. 
$\hfill\Box${\small
{\small

\end{document}